\documentclass[12pt]{article}

\usepackage{amssymb}
\usepackage{graphicx}
\usepackage{amsmath}
\usepackage{array}

\begin{document}

\newcommand{\B}{{\cal B}}
\newcommand{\C}{{\bf C}}
\newcommand{\SM}{{\bf SM}}
\newcommand{\T}{{\bf T}}

\newcommand{\F}{\mathbb{F}}
\newcommand{\Q}{\mathbb{Q}}

\newtheorem{theorem}{Theorem}[section]
\newtheorem{proposition}[theorem]{Proposition}
\newtheorem{lemma}[theorem]{Lemma}
\newtheorem{corollary}[theorem]{Corollary}
\newtheorem{problem}[theorem]{Problem}

\newcommand{\h}{\mathop{\rm ht }\nolimits}
\newcommand{\minh}{\mathop{\rm minht }\nolimits}
\newcommand{\maxh}{\mathop{\rm maxht }\nolimits}
\renewcommand{\r}{\mathop{\rm r }\nolimits}

\newcommand{\fillbx}{\rule[-0.1mm]{2.5mm}{2.5mm}}
\newcommand{\bbox}{\hfill\fillbx}
\newcommand{\emptybx}{\framebox[5mm]{\rule{3mm}{3mm}}}
\newcommand{\emptybox}{\hfill\emptybx}

\title{The Catalan matroid.}
\author{Federico Ardila \\ fardila@math.mit.edu}
\date{September 4, 2002}
\maketitle

\begin{abstract}
We show how the set of Dyck paths of length $2n$ naturally gives
rise to a matroid, which we call the ``Catalan matroid" $\C_n$.
We describe this matroid in detail; among several other results,
we show that $\C_n$ is self-dual, it is representable over $\Q$
but not over finite fields $\F_q$ with $q \leq n-2$, and it has a
nice Tutte polynomial.

We then generalize our construction to obtain a family of
matroids, which we call ``shifted matroids". They arose
independently and almost simultaneously in the work of Klivans,
who showed that they are precisely the matroids whose
independence complex is a shifted complex.

\end{abstract}

\section{Introduction}\label{sec:intro}

A Dyck path of length $2n$ is a path in the plane from $(0,0)$ to
$(2n,0)$, with steps $(1,1)$ and $(1,-1)$, that never passes
below the $x$-axis. It is a classical result (see for example
\cite[Corollary 6.2.3.(iv)]{St99}) that the number of Dyck paths
of length $2n$ is equal to the Catalan number $C_n =
\frac{1}{n+1} {2n \choose n}$.

Each Dyck path $P$ defines an \emph{up-step set}, consisting of
the integers $i$ for which the $i$-th step of $P$ is $(1,1)$. The
starting point of this paper is Theorem \ref{th:bases}. It states
that the collection of up-step sets of all Dyck paths of length
$2n$ is the collection of bases of a matroid. Most of this paper
is devoted to the study of this matroid, which we call the
\emph{Catalan matroid}, and denote $\C_n$.

Section \ref{sec:matroid} starts by proving Theorem
\ref{th:bases}. As we know, there are many equivalent ways of
defining a matroid: in terms of its rank function, its
independent sets, its flats, and its circuits, among others. The
rest of Section \ref{sec:matroid} is devoted to describing some
of these definitions for $\C_n$.

In Section \ref{sec:Tutte}, we compute the Tutte polynomial of
the Catalan matroid. We find that it enumerates Dyck paths
according to two simple statistics. Some nice enumerative results
are derived as a consequence.

In Section \ref{sec:shifted}, we generalize our construction of
$\C_n$ to a wider class of matroids, which we call \emph{shifted
matroids}. Their name is justified by a result of Klivans, who
discovered them independently, proving that they are precisely
those matroids whose independence complex is a shifted complex.
We then generalize our construction in a different direction to
obtain, for any finite poset $P$ and any order ideal $I$, a
shifted family of sets. This family is not always the set of
bases of a matroid.

Finally, in Section \ref{sec:rep} we address the question of
representability of the matroids we have constructed. We show
that the Catalan matroid, and more generally any shifted matroid,
is representable over $\Q$. In the opposite direction, we show
that $\C_n$ is \emph{not} representable over the finite field
$\F_q$ if $q \leq n-2$.

Throughout this paper, we will assume some familiarity with the
basic concepts of matroid theory. For instance, Chapter 1 of
\cite{Ox92} should be enough to understand most of the paper. We
also highly recommend Section 6.2 and exercises 6.19-6.37 of
\cite{St99} for an encyclopedic treatment of Catalan numbers
and related topics.

\section{The matroid}\label{sec:matroid}

Let $n$ be a fixed positive integer. Consider all paths in the
plane which start at the origin and consist of $2n$ steps, where
each step is either $(1,1)$ or $(1,-1)$. We will call such steps
\emph{up-steps} and \emph{down-steps}, respectively. From now on,
the word \emph{path} will always to refer to a path of this form.

Such paths are in bijection with subsets of $[2n]$. To each path
$P$, we can assign the set of integers $i$ for which the $i$-th
step of $P$ is an up-step. We call this set the \emph{up-step set}
of $P$. Conversely, to each subset $A \subseteq [2n]$, we can
assign the path whose $i$-th step is an up-step if and only if
$i$ is in $A$.

To simplify the notation later on, we will omit the brackets when
we talk about subsets of $[2n]$. We will also use subsets of
$[2n]$ and paths interchangeably. For example, for $n=3$, the path
$13$ will be the path with up-steps at steps $1$ and $3$, and
down-steps at steps $2, 4, 5$ and $6$.

A useful statistic to keep track of will be the \emph{height of
path $P$ at $x$}; \emph{i.e.}, the height of the path after
taking its first $x$ steps. We shall denote it $\h_P(x)$; it is
equal to $2|P_{\leq x}| - x$, where $P_{\leq x}$ denotes the set
of elements of $P$ which are less than or equal to $x$. Also, let
$\minh_P$ and $\maxh_P$ be the minimum and maximum heights that
$P$ achieves, respectively.

\begin{theorem} \label{th:bases}
Let $\B_n$ be the collection of up-step sets of all Dyck paths of
length $2n$. Then $\B_n$ is the collection of bases of a matroid.
\end{theorem}

\noindent \emph{Proof.} We need to check the two axioms for the
collection of bases of a matroid:

\begin{enumerate}
\item[{\bf (B1)}] $\B_n$ is non-empty.
\item[{\bf (B2)}] If $A$ and $B$ are members of $\B_n$ and
$a \in A - B$, then there is an element $b \in B - A$ such that
$(A - a) \cup b \in \B_n$.
\end{enumerate}

The first axiom is satisfied trivially, so we only need to check
the second one. Let $A$ and $B$ be members of $\B_n$, and let $a
\in A - B$. First we will describe those $k$ not in $A$ for which
$A \, - \, a \, \cup \, k \in \B_n$, and then we will show that
the smallest element of $B-A$ is one of them.

For $k \notin A$, consider the path $A \, - \, a \, \cup \, k$,
which is a very slight deformation of the path $A$. It still
consists of $n$ up-steps and $n$ down-steps; to determine if it
is a Dyck path, we just need to check whether it goes below the
$x$-axis. There are two cases to consider.

The first case is that $k<a$. In this case, for $k \leq c < a$, we
have that $\h_{A \, - \, a \, \cup \, k}(c) = \h_{A}(c) \, +  \,
2$. For all other values of $c$, we have that $\h_{A \, - \, a \,
\cup \, k}(c) = \h_{A}(c)$. Hence the path $A \, - \, a \, \cup
\, k$ stays above the path $A$, so it is Dyck.

The second case is that $a<k$. Here, for $a \leq c < k$, we have
that $\h_{A \, - \, a \, \cup \, k}(c) = \h_{A}(c) \, - \, 2$.
For all other values of $c$, we have that $\h_{A \, - \, a \,
\cup \, k}(c) = \h_{A}(c)$. Therefore, the path $A \, - \, a \,
\cup \, k$ is Dyck if and only if $\h_{A}(c) \geq 2$ for all $a
\leq c < k$.

With that simple analysis, we can show that $A \, - \, a \, \cup
\ b \in \B_n$, where $b$ is the smallest element of $B-A$. If $b
< a$, then we are done by the first case of our analysis.
Otherwise, consider an arbitrary $c$ with $a \leq c < b$. There
are no elements of $B-A$ less than or equal to $c$; so up to the
$c$-th step, every step which is an up-step in $B$ is also an
up-step in $A$. Furthermore, the $a$-th step is a down-step in
$B$ and an up-step in $A$. Therefore, $\h_A(c) \geq \h_B(c) \, +
\, 2 \geq 2$. This concludes our proof. \bbox

\medskip

A matroid is uniquely determined by its collection of bases.
Therefore Theorem \ref{th:bases} defines a matroid, which will
call \emph{the Catalan matroid of rank $n$} (or simply \emph{the
Catalan matroid}), and denote by $\C_n$. This paper is mostly
devoted to the study of this matroid.

\begin{proposition} \label{prop:rank}
The rank function of $\C_n$ is given by
$$
\r(A) = n + \lfloor \minh_A /2 \rfloor
$$
for each $A \subseteq [2n]$.
\end{proposition}

\noindent \emph{Proof.} Fix a subset $A \subseteq [2n]$, and let
$\minh_A = -y$, where $y$ is a non-negative integer. Also, let $x$
be the smallest integer such that $\h_A(x) = \minh_A$.

Recall that
the rank of a subset $A$ of $[2n]$ is equal to the largest
possible size of an intersection $A \cap B$, where $B$ is a basis
of $\C_n$.

The path $A$ is at height $-y$ after taking $|A_{\leq \,x}|$
up-steps and $x - |A_{\leq \, x}|$ down-steps, so $|A_{\leq \,x}|
= (x-y)/2$. Also, for any basis $B$, we have that $|B_{>\,x}| \leq
n-x/2$, since $\h_B(x) \geq 0$.  Hence
\begin{eqnarray*}
|A \cap B| &=& |(A \cap B)_{\leq \, x}| + |(A \cap B)_{>\,x}| \\
& \leq & |A_{\leq \, x}| + |B_{>\,x}| \leq n - y/2.
\end{eqnarray*}
We conclude that $\r(A) \leq n + \lfloor \minh_A /2 \rfloor$.

Now we need a basis $B$ with $|A \cap B| = n + \lfloor \minh_A /2
\rfloor$. We construct it as follows. First, add to $A$ the
smallest $a = \lceil y/2 \rceil$ numbers that it is missing, to
obtain the set $A'$. Then $\h_{A'}(x) = 2a - y \geq 0$; in fact,
it is clear that the path $A'$ never crosses the $x$-axis.
Let $|A| = n+h$ for some integer $h$; then $\h_A(2n) = 2h$ and
$\h_{A'}(2n) = 2h+2a$. Now
remove from $A'$ the largest $h+a$ numbers that it contains, to
obtain the set $B$. It is again easy to see that the path $B$
never crosses the $x$-axis, and ends at $(2n,0)$. So $B$ is Dyck, and
$$
|A \cap B| = |A \cap A'| - (h+a) = |A| - (h+a) = n-a
$$
as desired. \bbox

\medskip

Now that we know the rank function of $\C_n$, we can describe
several important classes of subsets of the matroid. We will only
provide a proof for the description of the class of flats; the
remaining proofs are similar in flavor. The interested
reader may want to complete the details to get better acquainted
with the matroid $\C_n$.

\begin{proposition} \label{prop:flats}
The flats of $\C_n$ are the subsets $A \subseteq [2n]$ such that
\begin{enumerate}
\item[(i)]
$\minh_A$ is odd, and
\item[(ii)]
if $\h_A(x) = \minh_A$, then $\{x+1, \ldots, 2n\} \subseteq A$.
\end{enumerate}
\end{proposition}

\noindent \emph{Proof.} Let $A$ be a flat of $\C_n$, and let $x$
be such that $\h_A(x) = \minh_A$. If some integer $y$ with $x+1
\leq y \leq n$ was not in $A$, then we would clearly have
$\minh_{A \, \cup y} = \minh_A$ and thus $\r(A \cup y) = \r(A)$,
contradicting the assumption that $A$ is a flat. Therefore, any
flat must satisfy condition \emph{(ii)}.

Also, if we had a flat $A$ with $\minh_A = -2h$ achieved at
$\h_A(x)$, then we would have $x \notin A$, and $\minh_{A \, \cup
x} = -2h+1$ would be achieved at $\h_{A \, \cup x}(x-1)$. We
would then have $\r(A \cup x) = \r(A)$, again a contradiction. So
any flat $A$ must also satisfy condition \emph{(i)}.

Conversely, assume that $A$ satisfies conditions \emph{(i)} and
\emph{(ii)}. Let $\minh_A = -(2k+1)$, which can only be achieved
once, say at $\h_A(x)$. Any $y$ which is not in $A$ must be less
than or equal to $x$; and we have $\minh_{A \cup y} = -(2k-1)$ if
$y<x$, or $\minh_{A \cup y} = -2k$ if $y=x$. In either case,
$\r(A \cup y) = \r(A)+1$. This completes the proof. \bbox

\medskip

\begin{proposition} \label{prop:ind}
The independent sets of $\C_n$ are the subsets $A \subseteq [2n]$
such that $\minh_A = \h_A(2n)$.
\end{proposition}

\begin{proposition} \label{prop:span}
The spanning sets of $\C_n$ are the subsets $A \subseteq [2n]$
such that $\minh_A = 0$.
\end{proposition}

\begin{proposition} \label{prop:circuits}
The circuits of $\C_n$ are the subsets $A \subseteq [2n]$ of the
form $A = \{2k, 2k + b_1, \ldots, 2k + b_{n-k}\}$, for some
positive integer $k \leq n$ and some Dyck path $\{b_1, \ldots,
b_{n-k}\}$ of length $2(n-k)$.
\end{proposition}

\begin{proposition} \label{prop:bonds}
The bonds of $\C_n$ are the subsets $A \subseteq [2n]$ such that
\begin{enumerate}
\item[(i)]
$\maxh_A=1$, and
\item[(ii)]
if $\h_A(x) = 1$, then $A$ has no elements greater than $x$.
\end{enumerate}

\end{proposition}

We complete this section with an observation which is interesting
in itself, and will also be important to us in section
\ref{sec:Tutte}.

\begin{proposition} \label{prop:dual}
The Catalan matroid is self-dual.\footnote{We follow Oxley
\cite{Ox92} in calling a matroid $M$ \emph{self-dual} if $M \cong
M^*$. It is worth mentioning, however, that some authors reserve
the term `self-dual' for matroids $M$ such that $M = M^*$.}
\end{proposition}

\noindent \emph{Proof.} Say $B = \{b_1, \ldots, b_n\}$ is a basis
of $\C_n$, and let $[2n]-B = \{c_1, \ldots, c_n\}$ be the
corresponding basis of the dual matroid $\C_n^{\,*}$. Then
$\{2n+1-c_n, \ldots, 2n+1-c_1\}$ is a Dyck path; in fact, it is
the path obtained by reflecting the Dyck path $B$ across a
vertical axis. So the bases of $\C_n^{\,*}$ are simply the up-step
sets of all Dyck paths of length $2n$, under the relabeling $x
\rightarrow 2n+1-x$. Thus $\C_n^{\,*} \cong \C_n$. \bbox

\medskip

\section{The Tutte polynomial} \label{sec:Tutte}

Given a matroid $M$ over a ground set $S$, its Tutte polynomial
is defined as:
$$
T_{M}(q,t) = \sum_{A \subseteq S}
(q-1)^{\r(S)-\r(A)}(t-1)^{|A|-\r(A)}.
$$

For our purposes, it is more convenient to define the Tutte
polynomial in terms of the internal and external activity of the
bases. We recall this definition now.

We first need to fix an arbitrary linear ordering of $S$.

For any basis $B$ and any element $e \notin B$, the set $B \cup
e$ contains a unique circuit. If $e$ is the smallest element of
that circuit with respect to our fixed linear order, then we say
that $e$ is \emph{externally active} with respect to $B$. The
number of externally active elements with respect to $B$ is
called the \emph{external activity} of $B$; we shall denote it by
$e(B)$.

Dually, for any basis $B$ and any element $i \in B$, the set $S-B
\cup i$ contains a unique bond. If $i$ is the smallest element of
that bond, then we say that $i$ is \emph{internally active} with
respect to $B$. The number of internally active elements with
respect to $B$ is called the \emph{internal activity} of $B$; we
shall denote it by $i(B)$.\footnote{The internally active
elements with respect to a basis $B$ of $M$ are precisely the
externally active elements with respect to the basis $S-B$ of the
dual matroid $M^*$. That is why we say that internal activity and
external activity are dual concepts. }

\begin{proposition}(Crapo, \cite{Cr69})\label{prop:Tutte}
For any matroid $M$ and any linear order of its ground set,
$$
T_{M}(q,t) = \sum_{B \, \mathrm{basis}} q^{i(B)} t^{e(B)}.
$$
\end{proposition}

We will use Proposition \ref{prop:Tutte} to study the Tutte
polynomial of the Catalan matroid. The first thing to do is to
fix a linear order of its ground set, $[2n]$. We will use the
most natural choice: $1<2<\cdots<2n$. Now we compute the internal
and external activity of each basis of $\C_n$.

\begin{lemma}\label{lemma:int}
The internal activity of a Dyck path $B$ is equal to the number of
up-steps that $B$ takes before its first down-step.
\end{lemma}

\noindent \emph{Proof.} Let $i \in B$. The path $[2n]-B$ never
goes above height $0$; the path $[2n]-B \cup i$ goes up to height
$2$. Let $j$ be the smallest integer such that $\h_{[2n]\, - \, B
\, \cup \, i}(j) = 1$. Clearly $j \geq i$.

Let $D$ be the unique bond of $\C_n$ which can be obtained by
deleting some elements of $[2n]-B \cup i$. We cannot delete any
element less than or equal to $j$, or else the resulting path
will not reach height $1$. We must delete any element larger than
$j$ by Proposition \ref{prop:bonds}. So $D = ([2n]-B)_{\leq j}$.

Therefore, $i$ is the smallest element of $D$ if and only if $B$
contains all of $1, 2, \ldots, i-1$. This completes the proof.
\bbox

\medskip

\begin{lemma}\label{lemma:ext}
The external activity of a Dyck path $B$ is equal to the number
of positive integers $x$ for which $\h_B(x)=0$.
\end{lemma}

\noindent \emph{Proof.} Let $e \notin B$. The path $B \cup e$ ends
at height $2$; let $2k-1$ be the largest integer such that $\h_{B
\cup e}(2k-1) = 1.$ Clearly $2k-1 < e$.

We start by showing that the unique circuit $C$ of $\C_n$
contained in $B \cup e$ is $(B \cup e)_{\geq 2k}$.

Since $C \subseteq B \cup e$, we have that $\h_C(2n) - \h_C(2k-1)
\leq \h_{B \cup e}(2n) - \h_{B \cup e}(2k-1) = 1$. Equality holds
if and only if every up-step of $B \cup e$ after the $(2k-1)$-th
is also an up-step of $C$; \emph{i.e.}, when $(B \cup e)_{\geq 2k}
= C_{\geq 2k}$.

But it is clear from Proposition \ref{prop:circuits} that
$\h_C(2n) - \minh_C = 1$, and that $\minh_C$ is only achieved at
$\h_C(\min C - 1)$. So the above inequality can only hold if $\min C
= 2k$. Thus $C = C_{\geq 2k} = (B \cup e)_{\geq 2k}$
as desired.

Now we know that $\min C = 2k$, so $e$ is externally active if and
only if $e=2k$. If $\h_B(e)=0$, this is clearly the case. On the
other hand, if $\h_B(e) \geq 1$, then $\h_{B \cup e}(e-1) =
\h_B(e-1) \geq 2$, so this is not the case. This completes the
proof. \bbox

\medskip

\begin{theorem} \label{th:Tutte}
For a Dyck path $P$, let $a(P)$ denote the number of up-steps
that $P$ takes before its first down-step, and let $b(P)$ denote
the number of positive integers $x$ for which $\h_P(x)=0$

Then the Tutte polynomial of the Catalan matroid $\C_n$ is equal
to
$$
T_{\C_n}(q,t) = \sum_{P  \, \mathrm{Dyck}} q^{a(P)} t^{b(P)},
$$
where the sum is over all Dyck paths of length $2n$.
\end{theorem}

\noindent \emph{Proof.} This follows immediately from Proposition
\ref{prop:Tutte} and Lemmas \ref{lemma:int} and \ref{lemma:ext}.
\bbox

\medskip

\begin{corollary} \label{cor:sym}
The polynomial
$$
\sum_{P  \, \mathrm{Dyck}} q^{a(P)} t^{b(P)},
$$
is symmetric in $q$ and $t$.
\end{corollary}

\noindent \emph{Proof.} It is well-known that, for any matroid
$M$, we have $T_{M^*}(q,t) = T_M(t,q)$. The result follows from
Proposition \ref{prop:dual} and Theorem \ref{th:Tutte}. \bbox

\medskip

It is a known fact that the statistics $a(P)$ and $b(P)$ are
equidistributed over the set of Dyck paths of length $2n$. The
number of paths with $a(P)=k$ and the number of paths with
$b(P)=k$ are both equal to $\frac{k}{2n-k}{2n-k \choose n}$. For
the first equality, see for example \cite{Va97}; for the second,
see \cite[equation (7)]{Kr70}.

Corollary \ref{cor:sym} was also discovered independently by James
Haglund \cite{Ha02}. It is not difficult to prove it directly; in
fact, it will be an immediate consequence of our next theorem.

\begin{theorem}
Let $C(x) = \frac{1 \, - \, \sqrt{1 \, - \, 4x}}2 = C_0 + C_1x +
C_2x^2 + \cdots$ be the generating function for the Catalan
numbers. Then
$$
\sum_{n \geq 0} T_{\C_n}(q,t)x^n =
\frac{1+(qt-q-t)xC(x)}{1-qtx+(qt-q-t)xC(x)}.
$$
\end{theorem}

\noindent \emph{Proof.} A Dyck path $P$ of length $2n \geq 2$ can
be decomposed uniquely in the standard way: it starts with an
up-step, then it follows a Dyck path $P_1$ of length $2r$, then
it takes a down-step, and it ends with a Dyck path $P_2$ of length
$2s$, for some non-negative integers $r, s$ with $r+s = n-1$. More
precisely, and necessarily more confusingly,
$$
P = \{1, 1+p_1, 1+p_2, \ldots, 1+p_r, 2r+2+q_1, 2r+2+q_2, \ldots,
2r+2+q_s\}
$$
for some Dyck paths $\{p_1, \ldots, p_r\}$ and $\{q_1, \ldots,
q_s\}$ with $r+s=n-1$.

It is clear that in this decomposition we have $a(P) = a(P_1)+1$
and $b(P) = b(P_2)+1$. Therefore
\begin{eqnarray*}
T_{\C_n}(q,t) &=& \sum_{r+s=n-1} \, \sum_{P_1 \, \in  \, \B_r}
\sum_{P_2 \, \in \, \B_s} q^{a(P_1) + 1} \, t^{b(P_2)+1} \\
&=& qt \sum_{r+s=n-1} T_{\C_r}(q,1) T_{\C_s}(1,t)
\end{eqnarray*}
for $n \geq 1$; so if we write $\T (q,t,x) = \sum_{n \geq 0}
T_{\C_n}(q,t)x^n$, we have
\begin{equation} \label{eq:T}
\T(q,t,x) = 1 + qtx \T(q,1,x) \T(1,t,x).
\end{equation}
Now observe that $\T(1,1,x) = C(x)$. Setting $q=1$ in
(\ref{eq:T}) gives a formula for $\T(1,t,x)$, and setting $t=1$
gives a formula for $\T(q,1,x)$. Substituting these two formulas
back into (\ref{eq:T}), we get the desired result. \bbox

\medskip

\section{Shifted matroids} \label{sec:shifted}

We now generalize our construction of $\C_n$ to a larger family of
matroids, which we call \emph{shifted matroids}. There is one
shifted matroid for each non-empty set $S = \{s_1 < \cdots <
s_n\}$ of positive integers, which we shall denote $\SM(s_1,
\ldots, s_n)$

\begin{theorem} \label{th:shiftedmatroid}
Let $S = \{s_1 < \cdots < s_n\}$ be a set of positive integers,
and let $\B_S$ be the collection of sets $\{a_1 < \cdots < a_n\}$
such that $a_1 \leq s_1, \ldots, a_n \leq s_n$. Then $\B_S$ is the
collection of bases of a matroid $\SM(s_1, \ldots, s_n)$.
\end{theorem}

\noindent \emph{Proof.} Once again, as in the proof of Theorem
\ref{th:bases}, axiom {\bf (B1)} is trivial, since $S \in B_S$. We
need to check axiom {\bf (B2)}. Let $A=\{a_1 < \cdots < a_n\}$ and
$B=\{b_1 < \cdots < b_n\}$ be in $B_S$, and let $a_x \in A-B$. We
claim that, if $b_y$ is the smallest element in $B-A$, then $A -
a_x \cup b_y \in \B_S$.

Let $i$ be the integer such that $a_i < b_y < a_{i+1}$. (If
$b_y<a_1$ then the claim is trivially true, since we are
replacing $a_x$ in $A$ with a number smaller than it. If
$b_y>a_n$ then set $i=n$.) We may assume that $i \geq x$; if that
was not the case, then we would have $b_y < a_{i+1} \leq a_x$, and
the claim would be trivial. We then have
$$
A-a_x \cup b_y = \{a_1 < \cdots < a_{x-1} < a_{x+1} < \cdots < a_i
< b_y < a_{i+1} < \cdots < a_n \}
$$
and we have $n$ inequalities to check.

The first $x-1$ and the last $n-i$ of these do not require any
extra work: we already know that $a_k \leq s_k$ for $1 \leq k \leq
x-1$ and for $i+1 \leq k \leq n$.

For each value of $k$ with $x \leq k \leq i-1$, we need to check
that $a_{k+1} \leq s_k$. If $k \geq y$, we have $a_{k+1} \leq a_i
< b_y \leq b_k \leq s_k$. Otherwise, if $k<y$, proceed as follows.
Since $b_y$ is the smallest element of $B$ which is not in $A$,
and $a_x$ is not in $B$, the numbers $b_1, \ldots, b_k$ must all
be somewhere in the list $a_1, \ldots, a_{x-1}, a_{x+1}, \ldots
a_i$. Therefore the $k$-th smallest number of this list,
$a_{k+1}$, must be less than or equal to $b_k$, which is less
than or equal to $s_k$.

Finally, we need to check that $b_y \leq s_i$. Since the numbers
$b_1, \ldots, b_{y-1}$ all appear in the list $a_1, \ldots,
a_{x-1}, a_{x+1}, \ldots, a_k, \ldots, a_i$, we have $y-1 \leq
i-1$. Therefore $b_y \leq s_y \leq s_i$. \bbox

\medskip

A path $\{a_1 < \cdots < a_n\}$ is Dyck if and only if, for each
$i$ with $1 \leq i \leq n$, the $i$-th up-step comes before the
$i$-th down-step; that is, if and only if $a_i \leq 2i+1$.
Therefore, the Catalan matroid $\C_n$ is exactly the shifted
matroid $\SM(1,3,5,\ldots, 2n-1)$, with an additional loop $2n$.

\medskip

Recall that an \emph{abstract simplicial complex} $\Delta$ on
$[n]$ is a family of subsets of $[n]$ (called \emph{faces}) such
that if $G \in \Delta$ and $F \subseteq G$, then $F \in \Delta$.
A simplicial complex $\Delta$ is \emph{shifted} if, for any face
$F \in \Delta$ and any pair of elements $i<j$ such that $i \notin
F$ and $j \in F$, the subset $F - j \cup i$ is also a face of
$\Delta$.

The family of independent sets of a matroid $M$ is always a
simplicial complex, called the \emph{independence complex} of
$M$. For shifted matroids, we have the following simple observation.

\begin{proposition} \label{prop:shifted}
The independence complex of a shifted matroid $\SM(s_1, \ldots,
s_n)$ is a shifted complex.
\end{proposition}

\noindent \emph{Proof.} If $F \subseteq [s_n]$ is independent, it
is contained in some basis $B$. Now assume that we have two
elements $i<j$ such that $i \notin F$ and $j \in F$, and let $G =
F - j \cup i$. If the basis $B$ contains $i$, then it contains
$G$. Otherwise, $B - j \cup i$ is also a basis: for any $1 \leq k
\leq n$, its $k$-th smallest is less than or equal to the $k$-th
smallest element of $B$, which is less than or equal to $s_k$.
This basis contains $G$. In both cases, we conclude that $G$ is
independent. \bbox

\medskip

In \cite{Kl02}, Klivans characterizes \emph{shifted matroid
complexes}: shifted complexes which are the independence complex
of a matroid. Her result and ours were discovered almost
simultaneously. When we sat down to discuss them, we realized
that the matroids that arise in her characterization are
precisely the ones in our construction. This is why they were
baptized ``shifted matroids".

\begin{proposition}(Klivans, \cite{Kl02}) \label{prop:Carly}
If the independence complex of a loop-less matroid $M$ is a
shifted complex, then $M \cong \SM(s_1, \ldots, s_n)$ for some
positive integers $s_1 < \cdots < s_n$.
\end{proposition}

\medskip

Theorem \ref{th:shiftedmatroid} and Propositions
\ref{prop:shifted} and \ref{prop:Carly} have a nice application to
Young tableaux. Recall that a \emph{partition} $\lambda =
(\lambda_1, \ldots, \lambda_k)$ of $n$ is a weakly decreasing
sequence of positive integers which add up to $n$. We associate
to it a \emph{Young diagram}: a left-justified array of unit
squares, which has $\lambda_i$ squares on the $i$-th row from top
to bottom.\footnote{This is the English way of drawing Young
diagrams; francophones draw them with $\lambda_i$ squares on the
$i$-th row from bottom to top.} A \emph{standard Young tableaux}
is a placement of the integers $1, \ldots, n$ in the squares of
the Young diagram, in such a way that the numbers are increasing
from left to right and from top to bottom.

These definitions will be sufficient for our purposes. For a much
deeper treatment of the theory of Young tableaux, we refer the
reader to \cite{Fu97}.

\begin{corollary} \label{cor:partition}
Let $\lambda$ be a partition. Define the \emph{first row set} of
a standard Young tableau $T$ of shape $\lambda$ to be the set of
entries which appear in the first row of $T$. Then the collection
of first row sets of all standard Young tableaux of shape
$\lambda$ is the collection of bases of a shifted matroid.
\end{corollary}

\noindent \emph{Proof.} Let $\lambda' = (\lambda'_1, \ldots,
\lambda'_n)$ be the conjugate partition of $\lambda$, so
$\lambda'_i$ is the number of squares on the $i$-th column of the
Young diagram of $\lambda$. Let $s_i = 1 + \lambda'_1+ \cdots +
\lambda'_{i-1}$ for $1 \leq i \leq n$.

Let $\{b_1 < \cdots < b_n\}$ be the first row set of a standard
Young tableau $T$ of shape $\lambda$. The first entry on the
$i$-th column of $T$ is $b_i$; it is smaller than every entry to
its southeast.  There are only $\lambda'_1+ \cdots +
\lambda'_{i-1}$ cells which are not to its southeast, so $b_i
\leq s_i$.

Conversely, if $B = \{b_1 < \cdots < b_n\}$ is such that $b_i \leq
s_i$ for $1 \leq i \leq n$, then we can construct a standard Young
tableau with first row set $B$. To do it, we first put the
elements of $B$ in order on the first row of $\lambda$. Then we
put the remaining numbers from $1$ to $|\lambda|$ on the
remaining cells going in order down the columns, starting with
the leftmost column. The inequalities $b_i \leq s_i$ guarantee
that this process does indeed give a Young tableau $T$.

It follows that the collection in question is simply the
collection of bases of the matroid $\SM(s_1, \ldots, s_n)$. \bbox

\medskip

We might try to generalize Corollary $\ref{cor:partition}$,
replacing the first row of $\lambda$ by any partition $\mu
\subseteq \lambda$. Define the \emph{$\mu$-set} of a standard
Young tableau $T$ of shape $\lambda$ to be the set of entries
which appear in the sub-shape $\mu$ in $T$.

It is not too difficult to see that we do not always get the
collection of bases of a matroid with this construction. However,
we can still say something interesting.

\begin{proposition} \label{prop:shiftedfamily}
Let $\mu \subseteq \lambda$ be partitions. Then the collection
$\B_{\lambda \mu}$ of $\mu$-sets of all standard Young tableau of
shape $\lambda$ is a shifted family.
\end{proposition}

\noindent \emph{Proof.} In fact, we prove something more general.
Let $P$ be a partially ordered set, or \emph{poset}, of $n$
elements. Recall that a subset $I$ of $P$ is an \emph{order
ideal} of $P$ if, for any pair of elements $x,y \in P$ with  $x
<_P y$ and $y \in I$, we also have $x \in I$. Also recall that a
\emph{linear extension} of $P$ is a bijection $f:P \rightarrow
[n]$ such that $i <_P j$ implies that $f(i)<f(j)$. For more
information on posets, we refer the reader to \cite[Chapter
3]{St86}.

Define the \emph{I-set} of a linear extension $f$ of $P$ to be
the set $\{f(i) : i \in I\}$.

\begin{proposition} \label{prop:P,I}
Let $P$ be a poset of $n$ elements, and let $I$ be an order ideal
of $P$ . Then the collection $\B_{P, \, I}$ of $I$-sets of all
linear extensions of $P$ is a shifted family.
\end{proposition}

\noindent \emph{Proof of Proposition \ref{prop:P,I}.} We need to
check that if we have a set $B \in \B_{P, \, I}$ and a pair of
numbers $a < b$ such that $a \notin B$ and $b \in B$, then $B - b
\cup a \in \B_{P, \, I}$. It is enough to show this for $a =
b-1$; the general case will then follow by induction on $b-a$.

So let $f$ be a linear extension of $P$ with $I$-set $B$, and let
$b \in B$ be such that $b-1 \notin B$. Let $b=f(i)$ and $b-1 =
f(p)$ where $i \in I$ and $p \in P-I$. Let $g:P \rightarrow [n]$
be defined by switching the values of $f$ at $i$ and $p$\,;
\emph{i.e.},
\begin{equation}
g(x) = \left\{ \begin{array}{ll}
                f(x) & \mbox{if $x \notin \{i,p\}$} \\
                b-1 & \mbox{if $x = i$} \\
                b & \mbox{if $x=p$}
                \end{array} \right.
\end{equation}

We claim that $g$ is also a linear extension for $P$. An important
observation is that $i$ and $p$ are incomparable in $P$. If we
had $i<p$, then we would have $b = f(i) < f(p) = b-1$. If we had
$i > p$, then $i \in I$ would imply $p \in I$.

To check that $f$ is a linear extension, we need to check that
$f(i) = b$ satisfies several inequalities: it must be greater than
all the values that $f$ takes on $P_{<\,i}$, and less than all
the values that $f$ takes on $P_{>\,i}$. But $b$ is never
compared to $b-1$ here, since $p$ and $i$ are incomparable.
Therefore, $b-1$ also satisfies all those inequalities that $b$
needs to satisfy.

Similarly, $b-1$ must be greater than all the values that $f$
takes on $P_{<p}$ and less than all the values that $f$ takes on
$P_{>p}$. The number $b$ also satisfies these inequalities.

So we can switch the values of $f(i)$ and $f(p)$, and the
resulting function $g$ will also be a linear extension of $P$.
Also, the $I$-set of $g$ is $B - b \cup (b-1)$. This concludes the
proof. \bbox

\medskip

Now, to prove Proposition \ref{prop:shiftedfamily}, partially
order the cells of $\lambda$: cell $i$ is less than cell $j$ in
$P_{\lambda}$ if and only if cell $i$ is northeast of cell $j$ in
$\lambda$. The cells of $\mu$ define an order ideal $I_{\mu}$ of
this poset $P_{\lambda}$, and $\B_{\lambda \mu} = \B_{P_{\lambda},
\, I_{\mu}}.$ Now use Proposition \ref{prop:P,I}. \bbox

\section{Representability} \label{sec:rep}

A natural question to ask is whether the Catalan matroid can be
represented as the vector matroid of a collection of vectors. We
answer that question in this section.

Given a collection of real numbers $x_1, \ldots, x_k$, let $x_S =
\prod_{i \in S} x_i$ for each subset $S \subseteq [k]$. Form all
the $2^{2^k}$ possible sums of some of the $x_S$'s. If these sums
are all distinct, we will say that the initial collection of
numbers is \emph{generic}. Most collections of real numbers are
generic. A specific example is a set of algebraically
independent real numbers. Another example is any sequence of
positive integers which increases quickly enough; for instance,
one that satisfies $x_i > (1+x_1)(1+x_2) \cdots (1+x_{i-1})$ for
$1 < i \leq k$.

\begin{theorem} \label{th:rep}
Let $v_1, \ldots, v_{2n}$ be the columns of a matrix
\[
A =
\begin{pmatrix}
a_{11} &   0    &    0   &    0   &    0   &   0    & \hdots &    0   &    0 \\
a_{21} & a_{22} & a_{23} &    0   &    0   &   0    & \hdots &    0   &    0 \\
a_{31} & a_{32} & a_{33} & a_{34} & a_{35} &   0    & \hdots &    0   &    0 \\
\vdots & \vdots & \vdots & \vdots & \vdots & \vdots & \ddots & \vdots & \vdots\\
a_{n1} & a_{n2} & a_{n3} & a_{n4} & a_{n5} & a_{n6} & \hdots & a_{n, \, 2n-1}&    0 \\
\end{pmatrix}
\]
where the $a_{ij}$'s with $1 \leq i \leq n$ and $1 \leq j \leq
2i-1$ are generic integers. Then the vector matroid of $\{v_1,
\ldots, v_n\}$ is isomorphic to the Catalan matroid $\C_n$.
\end{theorem}

\noindent \emph{Proof.} Let $M$ be the vector matroid of $V$. Let
$1 \leq b_1 < \cdots < b_n \leq 2n$. The set $B = \{v_{b_1},
\ldots, v_{b_n}\}$ is a basis for $M$ if and only if it is
independent; that is, if and only if the determinant of the
matrix $A_B$ with columns $v_{b_1}, \ldots, v_{b_n}$ is non-zero.

This determinant is a sum of $n!$ terms, with plus or minus signs
attached to them. Since the $a_{ij}$'s are generic, this sum can
only be zero if all the terms are $0$. So $B$ is a basis as long
as at least one of the $n!$ terms in this determinant is non-zero.

The question is now whether it is possible to place $n$
non-attacking rooks on the non-zero entries of $A_B$; that is, to
choose $n$ non-zero entries with no two on the same row or column.
The marriage theorem \cite[Theorem 5.1]{va92} would be the
standard tool to attack this kind of question. However, $A_B$ is
such that any entry below or to the left of a non-zero entry is
also non-zero. This fact will make our argument shorter and
self-contained.

If $b_i \leq 2i-1$ for all integers $i$ with $1 \leq i \leq n$,
then the $(i,i)$ entry of $A_B$ is $a_{i,b_i} \neq 0$. Therefore
we can place $n$ non-attacking rooks on non-zero entries of $A_B$
by putting them on the main diagonal.

Conversely, suppose that we have a placement of $n$ non-attacking
rooks on non-zero entries of $A_B$. Let $i$ be any integer
between $1$ and $n$. Then the rooks on the first $i$ rows must be
on $i$ different columns. From the shape that the non-zero entries
of $A_B$ form, we conclude that the $i$-th row must contain $i$
different non-zero entries. Thus the $(i,i)$ entry of $A_B$, which
is precisely $a_{i, b_i}$, must be non-zero. Therefore $b_i \leq
2i-1$. \bbox
\medskip

The above proof generalizes immediately to any shifted matroid
$\SM(s_1, \ldots, s_n)$.

\begin{theorem}
Let $s_1 < \cdots < s_n$ be arbitrary positive integers. Let $v_1,
\ldots, v_{s_n}$ be the columns of a matrix $A=(a_{ij})_{1 \leq i
\leq n \, , \, 1 \leq j \leq s_n}$, where the $a_{ij}$'s with $1
\leq i \leq n$ and $1 \leq j \leq s_i$ are generic, and the
remaining $a_{ij}$'s are equal to $0$. Then the vector matroid of
$\{v_1, \ldots, v_{s_n}\}$ is isomorphic to the shifted matroid
$\SM(s_1, \ldots, s_n)$.
\end{theorem}

\medskip

Theorem \ref{th:rep} shows that the Catalan matroid is
representable over $\Q$, or even over a sufficiently large finite
field. In the other direction, we now show a negative result about
representing $\C_n$ over finite fields.

\begin{proposition} \label{prop:repq}
The Catalan matroid $\C_n$ is \emph{not} representable over the
finite field $\F_q$ if $q \leq n-2$.
\end{proposition}

\noindent \emph{Proof.} It is known (\cite{Ox92}, Proposition
6.5.2) and easy to show that the uniform matroid $U_{2,k}$ is
$\F_q$-representable if and only if $q \geq k-1$. A matroid
containing it as a minor is not representable over $\F_q$ for
$q \leq k-2$.
This suggests that we should find the largest $k$ for which
$U_{2,k}$ is a minor of $\C_n$.

We can use the Scum theorem (Higgs, \cite{Ox92}, Proposition
3.3.7), which essentially says that, if a matroid has a certain
minor, then it must have that minor hanging from the top of its
lattice of flats. Our question is then equivalent to finding the
largest $k$ for which there exists a rank\,-\,$(n-2)$ flat which is
contained in $k$ rank\,-\,$(n-1)$ flats.

\begin{lemma} \label{lemma:flats}
Let $A$ be a rank\,-\,$(n-2)$ flat, and let $x$ be the smallest
integer such that $\h_A(x) = -1$. Then there are exactly
$\frac{x+3}2$ rank\,-\,$(n-1)$ flats containing $A$.
\end{lemma}

\noindent \emph{Proof of Lemma \ref{lemma:flats}.} We know from
Propositions \ref{prop:rank} and \ref{prop:flats} that $\minh_A =
-3$ and that, once the path $A$ reaches height $-3$, say at
$\h_A(y)$, it only takes up-steps. We want to add elements to $A$
to obtain a path which reaches a minimum height $-1$, and only
takes up-steps after that.

Say that we add one element $a$ to $A$. This new up-step at $a$
comes before the $y$-th, so $\h_{A \cup a}(y) = -1$. If we don't
want to add any more elements to $A$, we have to make sure that
$A \cup a$ only reaches height $-1$ at $y$. For this to be true,
we need the new up-step $a$ to occur on or before the $x$-th step.
In $A$, there are $\frac{x+1}2$ down-steps up to the $x$-th to
choose from. Each one of these gives a rank\,-\,$(n-1)$ flat
containing $A$

On the other hand, if we are to add more elements to $A$ to
obtain a rank-$(r-1)$ flat $B$, they will all be less than $y$ so
we will have $\h_B(y)>0$. The minimum height in $B$ must then be
achieved at some $z$ for which $\h_A(z) = -1$. In fact, for this
$z$ to be unique, it must be the leftmost one, \emph{i.e.}, it
must be $x$. So the only possibility is that $B = A_{\leq x} \cup
\{x+1, \ldots, 2n\}$, which is indeed a rank\,-\,$(n-1)$ flat.
This concludes the proof of Lemma \ref{lemma:flats}. \bbox

\medskip

Having shown Lemma \ref{lemma:flats}, the rest is easy. The
rank\,-\,$(n-2)$ flat which is contained in the largest number of
rank\,-\,$(n-1)$ flats, is the latest one to arrive to height
$-1$. This flat is clearly $\{1, 2, \ldots, n-3, n-2, 2n\}$, which
arrives to height $-1$ after $2n-3$ steps. It is contained in
exactly $n$ rank\,-\,$(n-1)$ flats.

Therefore $\C_n$ contains $U_{2,n}$ as a minor, and thus it is
not representable over a field $\F_q$ with $q \leq n-2$. \bbox

\medskip

\section{Acknowledgements}

I would like to thank Nantel Bergeron, James Haglund, Carly
Klivans and Richard Stanley for helpful discussions. In
particular, Richard Stanley pointed out that Proposition
\ref{prop:shiftedfamily} is true at the level of generality of
Proposition \ref{prop:P,I}.

\end{document}